\documentclass[amssymb,amstex,amsmath,11pt]{article}
\usepackage[margin=1.4in]{geometry}
\usepackage{amssymb,amsmath}              
\usepackage{amsmath,amssymb,latexsym,amsfonts,url}
\usepackage{graphicx}
\usepackage{amsthm}
\usepackage{enumitem}
\usepackage{amssymb}
\usepackage{color}
\usepackage{epstopdf}
\renewcommand{\epsilon}{\varepsilon}
\theoremstyle{plain}
\newtheorem{thm}{Theorem}[section]
\newtheorem{prop}[thm]{Proposition}
\newtheorem{cor}[thm]{Corollary}
\newtheorem{lem}[thm]{Lemma}

\newtheorem*{theorem*}{Theorem}
\newtheorem*{proposition*}{Proposition}
\theoremstyle{definition}

\theoremstyle{remark}

\newcommand{\dist}{{\rm dist}}

\newcommand{\trace}{{\rm tr}}

\newcommand{\tr}{{\rm{tr}}}

\newcommand{\hess}{{\rm Hess}}

\def\RR{{\mathbb R}}

\def\({\left(}
\def\){\right)}

\begin{document}
\begin{title}
{Free boundary constant mean curvature surfaces in a strictly convex three-manifold}
\end{title}
\begin{author}{Sung-Hong Min \and Keomkyo Seo}\end{author}

\date{\today}

\maketitle

\begin{abstract}
\noindent Let $C$ be a strictly convex domain in a $3$-dimensional Riemannian manifold with sectional curvature bounded above by a constant and let $\Sigma$ be a constant mean curvature surface with free boundary in $C$. We provide a pinching condition on the length of the traceless second fundamental form on $\Sigma$ which guarantees that the surface is homeomorphic to either a disk or an annulus. Furthermore, under the same pinching condition, we prove that if $C$ is a geodesic ball of $3$-dimensional space forms, then $\Sigma$ is either a spherical cap or a Delaunay surface. \\


\end{abstract}

\section{Introduction}

Let $C$ be a strictly convex domain in an $n$-dimensional Riemannian manifold $M$, i.e.,  the second fundamental form  of the boundary $\partial C$ is positive definite. A smooth compact hypersurface $\Sigma$ properly immersed in $C$ with $\partial \Sigma \subset \partial C$ is called a {\it free boundary cmc-H hypersurface} in $C$ if $\Sigma$ has constant mean curvature $H$ in $C$ and meets $\partial C$ orthogonally along its boundary $\partial \Sigma$. Free boundary cmc-$H$ hypersurfaces can be characterized variationally as critical points of the volume functional for volume-preserving variations of $\Sigma$ in $C$ with $\partial \Sigma \subset \partial C$. In particular, when $H=0$, such a hypersurface is called a {\it  free boundary minimal hypersurface}.

Since the study of free boundary minimal surfaces was initiated by Courant \cite{Courant} and Lewy \cite{Lewy} in the 1940s, this topic has been extensively investigated for a long time. We refer the readers to a nice survey \cite{Hildebrandt} by Hildebrandt for classical works on free boundary problems and a recent interesting survey \cite{Li} by Li for more recent results and some interesting questions about free boundary minimal surfaces in the unit ball. Various results on the existence of free boundary minimal disks in a convex domain of $\mathbb{R}^3$ were obtained by using geometric measure theory (see \cite{AS, GJ, Struwe} for example). Later, Fraser \cite{Fraser} showed a general existence of free boundary minimal disks in a Riemannian manifold. Regarding uniqueness of free boundary minimal surfaces, Nitsche \cite{Nitsche} obtained the following remarkable result.
\begin{theorem*}[Nitsche \cite{Nitsche}]
The only free boundary minimal disks in the unit ball $B^3 \subset \mathbb{R}^3$ are equatorial disks.
\end{theorem*}
This result was extended to free boundary minimal disks in geodesic balls in the $3$-dimensional space forms of constant curvature by Ros-Souam \cite{RS} and Souam \cite{Souam}. Recently, Fraser-Schoen \cite{FS} generalized Nitsche's result to free boundary minimal disks in geodesic balls of arbitrary dimension in space forms. In \cite{WW}, a higher dimensional analogue for free boundary minimal hypersurfaces in the $n$-dimensional Euclidean unit ball $B^n$ was obtained under some graphical conditions.

The next simplest example of free boundary minimal surface in the $3$-dimensional Euclidean unit ball $B^3$ is the critical catenoid, which is a piece of a catenoid in $\mathbb{R}^3$ intersecting $\partial B^3$ orthogonally. A number of characterizations of the critical catenoid have been obtained by numerous geometers. For example, refer to \cite{Devyver, SZ, Tran} for the Morse index estimate, \cite{AN, BV, CMV, FS2016, FS2020, KM, LX,  McGrath} for uniqueness results of the critical catenoid, and \cite{BB} for a variational characterization in terms of $2$-dimensional Hausdorff measure. We should mention that examples of free boundary minimal surfaces in $B^3$ are very rich. In fact, many techniques have been employed to construct new examples of free boundary minimal surfaces in $B^3$ (see \cite{CFS, FPZ, FS, FS2016, KL, KW, Ketover} for instance). See also \cite{LZ, MNS} for more general existence results for free boundary minimal surfaces in a Riemannian manifold.

On the other hand, the following famous gap theorem gives a characterization of the equator and the Clifford minimal hypersurface $C_{m,n}:= \mathbb{S}^m \left(\sqrt{\frac{m}{n}}\right)\times \mathbb{S}^{n-m}\left(\sqrt{\frac{n-m}{n}}\right)$
in $\mathbb{S}^{n+1}$ for $1\leq m \leq n-1$.

\begin{thm}[Chern-do Carmo-Kobayashi \cite{CDK}, Lawson \cite{Lawson}, Simons \cite{Simons}] \label{thm:Simons}
Let $\Sigma$ be a closed minimal hypersurface in the unit sphere $\mathbb{S}^{n+1}$. Assume that the second fundamental form $A$ on $\Sigma$ satisfies
$$|A|^2 \leq n.$$
Then
\begin{enumerate}[label=(\roman*)]
  \item either $|A|^2=0$ and $\Sigma$ is an equator;
  \item or $|A|^2=n$ and $\Sigma$ is a Clifford minimal hypersurface $C_{m,n}$ for $1\leq m \leq n-1$.
\end{enumerate}
\end{thm}
The pioneering work \cite{FS2011, FS2016} by Fraser and Schoen revealed that there are several similarities between free boundary minimal surfaces in a Euclidean unit ball $B^n$ and closed minimal surfaces in $\mathbb{S}^n$. From this perspective, Ambrozio-Nunes \cite{AN} obtained the following interesting result for free boundary minimal surfaces in the $3$-dimensional Euclidean unit ball $B^3$, which is similar to Theorem \ref{thm:Simons}.
\begin{thm}[Ambrozio-Nunes \cite{AN}] \label{thm: AN}
Let $\Sigma$ be a compact free boundary minimal surface in the $3$-dimensional Euclidean unit ball $B^3$. Assume that for every point $x\in \Sigma$,
$$|A|^2(x) \langle x, \nu(x) \rangle^2  \leq 2,$$
where $\nu(x)$ denotes the unit normal vector at the point $x$ and $A$ denotes the second fundamental form of $\Sigma$. Then
\begin{enumerate}[label=(\roman*)]
  \item either $|A|^2(x) \langle x, \nu(x) \rangle^2 =0$ and $\Sigma$ is an equatorial flat disk;
  \item or $|A|^2(x) \langle x, \nu(x) \rangle^2 =2$ at some point $p\in \Sigma$ and $\Sigma$ is a critical catenoid.
\end{enumerate}
\end{thm}

Very recently, Theorem \ref{thm: AN} has been extended in various directions. Li-Xiong \cite{LX} proved a similar gap result for free boundary minimal surfaces in geodesic balls of three-dimensional hyperbolic space and hemisphere. Barbosa-Viana \cite{BV} generalized Theorem \ref{thm: AN} to free boundary minimal surfaces of higher codimension in the $n$-dimensional Eucldiean unit ball $B^n$. Cavalcante-Mendes-Vit\'{o}rio \cite{CMV} obtained a topological gap result for compact free boundary submanifolds in the $n$-dimensional Euclidean unit ball $B^n$ which are not necessarily minimal. In \cite{BCP}, Barbosa-Cavalcante-Pereira proved an analogue of Theorem \ref{thm: AN} for free boundary constant mean curvature surfaces in the $3$-dimensional Euclidean unit ball $B^3$.

In this paper, we deal with some rigidity results about free boundary cmc-$H$ surfaces in a strictly convex domain of a $3$-dimensional Riemannian manifold. In Section $3$, we investigate the topology of such surfaces under a pinching condition on the length of the second fundamental form. To do this, we use the Hessian of the distance function in the ambient space to provide a pinching condition for a topological characterization of free boundary cmc-$H$ surfaces. More precisely, we obtain the following.

\begin{theorem*}[see Theorem \ref{thm:topology}]
Let $M^3$ be a $3$-dimensional Riemannian manifold with  sectional curvature $K_M \leq c$ for some constant $c$. Let $\Sigma$ be a free boundary cmc-H surface in a strictly convex domain $C\subset M^3$. Suppose that umbilic points on $\Sigma$ are isolated unless $\Sigma$ is totally umbilical. For the function $f$ defined as in (\ref{eqn:definition}) and the distance $r$ measured from a fixed point $p\in C$, assume that every point in $\Sigma$ satisfies
\begin{align*}
\frac{1}{2} |\Phi|^2 \left(\nabla_\nu f \right)^2 \leq \left(f'' +H  \nabla_\nu f \right)^2,
\end{align*}
where $\nabla_\nu f = \langle \nu, \nabla f \rangle$. Here $\nu$ and $\nabla$ denote the unit normal vector field on $\Sigma$ and the Levi-Civita connection on $M$, respectively. When $c>0$, assume further that $r<\frac{\pi}{2 \sqrt{c}}$. Then $\Sigma$ is homeomorphic to a disk or an annulus.
\end{theorem*}

\noindent Here we require the assumption that umbilic points of $\Sigma$ are isolated unless $\Sigma$ is totally umbilical. One reason is that Nitsche-type result for free boundary constant mean curvature disks is absent in geodesic balls of general ambient spaces, which is different from the result of Ambrozio-Nunes in the $3$-dimensional Euclidean unit ball. It would be interesting to classify cmc-$H$ surfaces which satisfy the assumption. In particular, in case where $C$ is a geodesic ball of $3$-dimensional space forms, it is well-known that the assumption is automatically satisfied by the holomorphicity of the usual Hopf differential on $\Sigma$. Using this, we are able to prove the following.

\begin{theorem*} [see Theorem \ref{thm:main}]
Let $\overline{M}^3(c)$ be a $3$-dimensional space form of constant curvature $c$. Let $\Sigma$ be a free boundary cmc-$H$ surface in the geodesic ball $B_R(p) \subset \overline{M}^3(c)$ of radius $R$ centered at $p$. When $c>0$, assume that $R<\frac{\pi}{2 \sqrt{c}}$.  For the function $f$ defined as in (\ref{eqn:definition}) and the distance $r$ measured from $p$, assume that every point in $\Sigma$ satisfies
\begin{align*}
\frac{1}{2} |\Phi|^2 \left(\nabla_\nu f \right)^2 \leq \left(f'' +H  \nabla_\nu f \right)^2,
\end{align*}
where $\nabla_\nu f = \langle \nu, \nabla f \rangle$. Here $\nu$ and $\nabla$ denote the unit normal vector field on $\Sigma$ and the Levi-Civita connection on $\overline{M}^3(c)$, respectively.  Then $\Sigma$ is either a spherical cap or a Delaunay surface.
\end{theorem*}

We note that, when $c=0$, this theorem is the same as the result of Ambrozio-Nunes \cite{AN} and Barbosa-Cavalcante-Pereira \cite{BCP}. In order to prove this theorem, we give a characterization of Delaunay surfaces as follows. Given a constant $H$, if $\gamma$ is a geodesic circle lying on a totally geodesic plane $P$ in a $3$-dimensional space form, then there exists a unique cmc-$H$ surface meeting $P$ orthogonally along $\gamma$, which is a Delaunay surface (see Lemma \ref{lem:rotation}).

\section{Preliminaries}

Let $\Sigma$ be an $n$-dimensional hypersurface in an $(n+1)$-dimensional Riemannian manifold $M$. Denote by $\nabla$ the Levi-Civita connection on $M$. Let $A$ be the second fundamental form of $\Sigma$. For all tangent vector fields $X, Y$ on $\Sigma$,
$$A(X, Y) = \langle \nabla_X Y, \nu \rangle,$$
where $\nu$ denotes the unit normal vector field on $\Sigma$. We recall that the {\it (normalized) mean curvature} $H$ of $\Sigma$ is defined by
$$ H=\frac{1}{n}\trace(A)=\frac{1}{n} \sum_{i=1}^n \lambda_i,$$
where each $\lambda_i$ is the principal curvature of $\Sigma$ for $1 \leq i \leq n$. The {\it traceless} second fundamental form $\Phi$ is defined by
$$\Phi=A-H\cdot g_{\scriptscriptstyle\Sigma},$$
where $g_{\scriptscriptstyle\Sigma}$ denotes the induced metric on $\Sigma$. The squared norm of the traceless second fundamental form can be computed in terms of that of the second fundamental form as follows:
\begin{align} \label{eqn:squred norm}
|\Phi|^2 = \displaystyle\sum_{i=1}^n \left( \lambda_i -H \right)^2 = |A|^2 -n H^2.
\end{align}
We say that a point $p \in \Sigma$ is {\it umbilical} if $\Phi(p)=0$. Note that all the principal curvatures on $\Sigma$ are equal at the umbilical points.

In the study of a Riemannian manifold of variable curvature, it is very helpful to use the Hessian comparison theorem for the distance, which can be stated as follows (see \cite{GW, SY} for example):
\begin{thm}[Hessian comparison theorem] \label{thm: classical Hessian comparison}
Let $M$ be a Riemannian manifold with sectional curvature $K_M$ bounded from above by some constant $c$ and $\overline{M}(c)$ be a space form of constant curvature $c$. Let $r(\cdot)=\dist(p, \cdot)$ and $\overline{r}(\cdot)=\dist(\overline{p}, \cdot)$ be the distance functions on $M$ and $\overline{M}(c)$ for fixed points $p\in M$ and $ \overline{p} \in \overline{M}(c)$, respectively. Assume that $\gamma$ is a geodesic from $p$ to $q$ and $v\in T_qM$ is a vector at $q$ perpendicular to $\gamma$. Denote by $\hess_{\scriptscriptstyle M}$ and $\hess_{\scriptscriptstyle\overline{M}}$ the Hessians on $M$ and $\overline{M}(c)$, respectively. Then
\begin{align} \label{eqn:classical hessian comparison}
\hess_{\scriptscriptstyle M} r (v,v) \geq \hess_{\scriptscriptstyle\overline{M}} \overline{r} (\overline{v},\overline{v}),
\end{align}
where $\overline{v}\in T_{\overline{q}}\overline{M}(c)$ is a vector at $\overline{q} \in \overline{M}(c)$ with $|v|=|\overline{v}|$ and $\overline{r}(\overline{q})=r(q)$, which is perpendicular to the geodesic $\overline{\gamma}\subset \overline{M}(c)$ from $\overline{p}$ to $\overline{q}$.
\end{thm}

Let $M^n$ be an $n$-dimensional Riemannian manifold with sectional curvature $K_M \leq c$ for some constant $c$ and let $r(x)$ be the distance from a fixed point $p$ to $x\in M^n$. Denote by $g$ and $\nabla$ the metric and the Levi-Civita connection on $M^n$, respectively. Given a smooth function $f: [0, \infty) \rightarrow \RR$, the Hessian of the function $f(r)$ on $M^n$ can be computed as follows:
\begin{align*}
\hess_{\scriptscriptstyle M} f(r) (X,Y) &= \nabla_{\scriptscriptstyle X} \nabla_{\scriptscriptstyle Y} f(r) - \nabla_{\scriptscriptstyle \nabla_X Y} f(r)\\
&=  \nabla_{\scriptscriptstyle X} \left( f'(r) \nabla_{\scriptscriptstyle Y} r \right) - f'(r) \nabla_{\scriptscriptstyle \nabla_X Y} r\\
&= f'(r) \nabla_{\scriptscriptstyle X} \nabla_{\scriptscriptstyle Y} r +  f''(r) \nabla_{\scriptscriptstyle X} r \nabla_{\scriptscriptstyle Y} r  - f'(r) \nabla_{\scriptscriptstyle \nabla_X Y} r\\
&= f'(r) \hess_{\scriptscriptstyle M} r (X,Y) + f''(r) \, \nabla r \otimes \nabla r (X,Y)
\end{align*}
for all tangent vector fields $X, Y$ on $M^n$. Thus
\begin{equation} \label{eqn:hess1}
\hess_{\scriptscriptstyle M} f(r) = f'(r) \hess_{\scriptscriptstyle M} r + f''(r) \, \nabla r \otimes \nabla r.
\end{equation}
In the same manner, the Hessian of $f(\overline{r})$ on an $n$-dimensional space form $\overline{M}^n(c)$ of constant curvature $c$ is given by
\begin{equation} \label{eqn:hess2}
\hess_{\scriptscriptstyle \overline{M}} f(\overline{r}) = f'(\overline{r}) \hess_{\scriptscriptstyle \overline{M}} \overline{r} + f''(\overline{r}) \, \overline{\nabla} \overline{r} \otimes \overline{\nabla} \overline{r},
\end{equation}
where $\overline{r}$ denotes the distance from a fixed point $\overline{p}$ in $\overline{M}^n(c)$ and $\overline{\nabla}$ denotes the Levi-Civita connection on $\overline{M}^n(c)$.

Assume that $\gamma \subset M$ is a geodesic from $p$ to $q$ and $\overline{\gamma} \subset \overline{M}(c)$ is a geodesic from $\overline{p}$ to $\overline{q}$ for $\overline{r}(\overline{q})=r(q)$. Consider an orthogonal splitting
$$T_q M= T_q \gamma \oplus N_q \gamma$$
into the tangent and normal spaces of $\gamma$ at $q$, respectively. Similarly, we have
$$T_{\overline{q}} \overline{M}(c)=T_{\overline{q}}\overline{\gamma} \oplus N_{\overline{q}}\overline{\gamma}.  $$
With respect to this splitting, we decompose any vectors $v\in  T_q M$ and $\overline{v} \in T_{\overline{q}} \overline{M}(c)$ as
$$v= v^T + v^N ~~{\rm and}~~ \overline{v}=\overline{v}^T + \overline{v}^N.$$
Using these notations, for $v\in T_q M$ and $\overline{v} \in T_{\overline{q}} \overline{M}(c)$ satisfying that $|v|=|\overline{v}|$ and $|v^N|=|\overline{v}^N|$ (thus $|v^T|= |\overline{v}^T|$), we have the following:
\begin{equation*}
\begin{cases}
 \hess_{\scriptscriptstyle M} f(r)(v,v) \geq \hess_{\scriptscriptstyle \overline{M}} f(\overline{r})(\overline{v}, \overline{v})&~~~~{\text{if }} f' \geq 0,\\
 \hess_{\scriptscriptstyle M} f(r)(v,v) \leq \hess_{\scriptscriptstyle \overline{M}} f(\overline{r})(\overline{v}, \overline{v}) &~~~~ {\text{if }} f' \leq 0.
\end{cases}
\end{equation*}
To see this, if $f'(r)\geq  0$ (resp. $\leq 0$), applying (\ref{eqn:classical hessian comparison}), (\ref{eqn:hess1}) and (\ref{eqn:hess2}) gives
\begin{align*}
\hess_{\scriptscriptstyle M} f(r) (v,v) &= \hess_{\scriptscriptstyle M} f(r) (v^T + v^N, v^T + v^N) \\
&= f'(r) \hess_{\scriptscriptstyle M} r (v^N, v^N) + f''(r) \, \nabla r \otimes \nabla r (v^T, v^T)\\
&\geq (\text{resp.}\leq)\, f'(r) \hess_{\scriptscriptstyle\overline{M}} \overline{r} (\overline{v}^N, \overline{v}^N) + f''(\overline{r}) \, |v^T|^2\\
&= f'(\overline{r}) \hess_{\scriptscriptstyle\overline{M}} \overline{r} (\overline{v}^N, \overline{v}^N) + f''(\overline{r}) \, |\overline{v}^T|^2 \\
&= \hess_{\scriptscriptstyle\overline{M}} f (\overline{r}) (\overline{v}, \overline{v}),
\end{align*}
where we used the fact that $\hess_{\scriptscriptstyle M} r(v^T,v^T)=0$ and $\hess_{\scriptscriptstyle M} r(v^T,v^N)=0$ in the second equality. For later use, define the function $f(r)$ on $M^n$ as follows:
\begin{equation} \label{eqn:definition}
f(r)=
\begin{cases}
 \frac{1}{2}r^2 &~~ {\text{if }} K_M \leq c=0,\\
 \cosh (k r) &~~{\text{if }} K_M \leq c=-k^2 {\text{ for }} k>0 ,\\
 \cos (k r) &~~{\text{if }} K_M \leq c=k^2  {\text{ and }} r < \frac{\pi}{2\sqrt{c}} {\text{ for }} k>0.\\
\end{cases}
\end{equation}
In the same way, we define the function $f(\overline{r})$ on the space form $\overline{M}^n(c)$. It is well-known that the Hessian of $f(\overline{r})$ on $\overline{M}^n(c)$ satisfies (see \cite{Petersen} for example)
$$\hess_{\scriptscriptstyle \overline{M}} f(\overline{r}) = f''(\overline{r}) \cdot \overline{g},$$
where $\overline{g}$ denotes the metric on $\overline{M}^n(c)$. In other words,
\begin{displaymath}
\left\{
\begin{array}{llll}
\hess_{\scriptscriptstyle \overline{M}} \frac{1}{2} \overline{r}^2&=& \overline{g} &~~\text{if }c=0,\\
\hess_{\scriptscriptstyle \overline{M}} \cosh k\overline{r}&=& k^2\cosh kr \cdot \overline{g} &~~\text{if } c=-k^2,\\
\hess_{\scriptscriptstyle \overline{M}} \cos k\overline{r}&=& -k^2\cos kr \cdot \overline{g} &~~\text{if } c=k^2.
\end{array}
\right.\end{displaymath}
Keeping these observations in mind, we can derive the following.
\begin{lem}
Let $(M, g)$ be a Riemannian manifold with sectional curvature $K_M$ bounded from above by some constant $c$. Let $r(\cdot)=\dist(p, \cdot)$  be the distance function on $M$ for a fixed point $p\in M$ and let the function $f$ be defined as in (\ref{eqn:definition}). Then we have at $q\in M$,
\begin{equation} \label{eqn:hess3}
\begin{cases}
\hess_{\scriptscriptstyle M} f(r) (v, v) \geq f''(r) \cdot g (v, v) &~~{\text{if }} c \leq 0,\\
\hess_{\scriptscriptstyle M} f(r) (v, v) \leq f''(r) \cdot g (v, v) &~~{\text{if }}  c>0 {\text{ and }} r < \frac{\pi}{2\sqrt{c}},
\end{cases}
\end{equation}
for all $v\in T_q M$.
\end{lem}
\begin{proof}
Consider a geodesic $\gamma \subset M$ from $p$ to $q$. As before, for a space form $\overline{M}(c)$ of constant curvature $c$ with  the distance $\overline{r}$ from a fixed point $\overline{p}$, let $\overline{\gamma} \subset \overline{M}(c)$ be a geodesic from $\overline{p}$ to $\overline{q}$ with $\overline{r}(\overline{q})=r(q)$.

Suppose that $c\leq 0$. Then for $v\in T_q M$ and $\overline{v} \in T_{\overline{q}} \overline{M}(c)$ satisfying that $|v|=|\overline{v}|$ and $|v^N|=|\overline{v}^N|$,
\begin{align*}
\hess_{\scriptscriptstyle M} f(r) (v,v) &\geq \hess_{\scriptscriptstyle\overline{M}} f (\overline{r}) (\overline{v}, \overline{v}) \\
&= f''(\overline{r}) \cdot \overline{g} (\overline{v}, \overline{v})\\
&= f''(r) \cdot g (v, v)
\end{align*}
at $q \in M$. One can similarly prove the case where $c>0$.
\end{proof}

\section{Topology of free boundary cmc-$H$ surfaces}
In this section, we provide a topological result for free boundary cmc-$H$ surfaces $\Sigma$ inside a strictly convex domain in a $3$-dimensional Riemannian manifold under the assumption that umbilic points on $\Sigma$ are isolated unless $\Sigma$ is totally umbilical and a pinching condition on $\Sigma$.

As mentioned in the introduction, Ambrozio-Nunes \cite{AN} obtained a gap theorem for free boundary minimal surfaces in a $3$-dimensional Euclidean unit ball $B^3$ with a certain pinching condition on the length of the second fundamental form. Note that the argument in \cite{AN} works also for properly embedded complete minimal surfaces  in $\mathbb{R}^3$. It turned out that any properly embedded complete minimal surface satisfying the same geometric condition must be either the plane or the catenoid, which was due to Meeks-P\'{e}rez-Ros \cite{MPR}.  Recently, there have been many interesting extensions \cite{BCP, BV, CMV, LX} of the work by Ambrozio-Nunes.  We extend the previous results into free boundary cmc-$H$ surfaces $\Sigma$ inside a strictly convex three-manifold under a similar pinching condition in terms of the distance function. Adopting the arguments of \cite{AN, BCP}, we prove the following.

\begin{thm} \label{thm:topology}
Let $M^3$ be a $3$-dimensional Riemannian manifold with  sectional curvature $K_M \leq c$ for some constant $c$. Let $\Sigma$ be a free boundary cmc-H surface in a strictly convex domain $C\subset M^3$. Suppose that umbilic points on $\Sigma$ are isolated unless $\Sigma$ is totally umbilical. For the function $f$ defined as in (\ref{eqn:definition}) and the distance $r$ measured from a fixed point $p\in C$, assume that every point in $\Sigma$ satisfies
\begin{align}
\frac{1}{2} |\Phi|^2 \left(\nabla_\nu f \right)^2 \leq \left(f'' +H  \nabla_\nu f \right)^2, \label{ineq:Assumption}
\end{align}
where $\nabla_\nu f = \langle \nu, \nabla f \rangle$. Here $\nu$ and $\nabla$ denote the unit normal vector field on $\Sigma$ and the Levi-Civita connection on $M$, respectively. When $c>0$, assume further that $r<\frac{\pi}{2 \sqrt{c}}$. Then $\Sigma$ is homeomorphic to either a disk or an annulus.
\end{thm}

\begin{proof}
Choose a local orthonormal frame $\{e_1, e_2, \nu \}$ on $M$ such that, restricted to $\Sigma$, the vectors $e_1, e_2$ are tangent to $\Sigma$ and the remaining vector $\nu$ is normal to $\Sigma$. Then the Hessian of the function $f$ on $\Sigma$ is given by
\begin{align} \label{eqn:hess4}
\hess_{\scriptscriptstyle\Sigma} f(e_i, e_j) &= e_i e_j f -\left( \nabla^{\scriptscriptstyle \Sigma}_{e_i} e_j \right) f \nonumber\\
&= e_i e_j f + \left( \nabla_{e_i} e_j - \nabla^{\scriptscriptstyle \Sigma}_{e_i} e_j \right) f -\left( \nabla_{e_i} e_j \right) f \nonumber\\
&= \left(e_i e_j  - \nabla_{e_i} e_j \right) f + \left(\nabla_{e_i} e_j  - \nabla^{\scriptscriptstyle \Sigma}_{e_i} e_j \right) f \nonumber\\
&= \hess_{\scriptscriptstyle M} f(e_i, e_j) +\left(\nabla_{e_i} e_j \right)^{\perp} f \nonumber\\
&= \hess_{\scriptscriptstyle M} f(e_i, e_j) +A(e_i, e_j) \nabla_\nu f,
\end{align}
where $\nabla^{\scriptscriptstyle \Sigma}$ denotes the induced connection on $\Sigma$ and $A$ denotes the second fundamental form of $\Sigma$. In order to prove Theorem \ref{thm:topology}, we need the following steps.\\

\noindent Step I: We claim that \emph{the geodesic curvature $\kappa_g$ of $\partial \Sigma$ is positive}. To see this, let $\alpha(s)$ be a parametrization of $\partial \Sigma$ by arclength parameter $s$. The assumption that $C$ is strictly convex implies that the curvature vector $\nabla_{\frac{d\alpha}{ds}} \frac{d\alpha}{ds}$ points inward. The geodesic curvature $\kappa_g$ of $\partial \Sigma$ in $\Sigma$ satisfies
$$\kappa_g = \langle \nabla_{\frac{d\alpha}{ds}} \textstyle\frac{d\alpha}{ds}, \eta \rangle>0,$$
where $\eta$ denotes the unit inward pointing conormal vector field along $\partial \Sigma$. In particular,  $\partial \Sigma$ is strictly convex in $\Sigma$.\\

\noindent Step II: We claim that  \emph{if either $\Sigma$ is totally umbilical or $\Sigma$ has nonnegative Gaussian curvature everywhere, then $\Sigma$ is homeomorphic to a disk}. To see this, consider the case where $\Sigma$ is totally umbilical. Then the Gaussian curvature $K_\Sigma$ of $\Sigma$ satisfies
$$K_\Sigma = H^2 \geq 0.$$
This implies that, in any case, $\Sigma$ has nonnegative Gaussian curvature everywhere. From the Gauss-Bonnet theorem and Step I, it follows
\begin{displaymath}
\int_{\Sigma} K_{\scriptscriptstyle \Sigma} + \int_{\partial \Sigma} \kappa_g = 2\pi \chi (\Sigma) >0,
\end{displaymath}
which shows that $\chi(\Sigma)=1$. Thus $\Sigma$ is orientable and has exactly one boundary component. Therefore we see that $\Sigma$ is homeomorphic to a disk.\\

By Step II, from now on, we may assume that $\Sigma$ is not a totally umbilical surface and has negative Gaussian curvature at some point of $\Sigma$.\\

\noindent Step III: We claim that, for $i=1, 2$,
\begin{equation*}
\begin{cases}
\hess_{\scriptscriptstyle \Sigma} f (e_i, e_i) \leq 0  &~~{\text{if }} c > 0,\\
\hess_{\scriptscriptstyle \Sigma} f (e_i, e_i) \geq 0 &~~{\text{if }}  c \leq 0.
\end{cases}
\end{equation*}
To see this, we first consider the case where $c>0$. Combining (\ref{eqn:hess3}) with (\ref{eqn:hess4}) gives
\begin{eqnarray*}
\hess_{\scriptscriptstyle\Sigma} f (e_i, e_i) &=& \hess_{\scriptscriptstyle M} f (e_i, e_i)+ A(e_i , e_i) \nabla_\nu f\\
&\leq& f'' + A(e_i , e_i) \nabla_\nu f
\end{eqnarray*}
for $i=1, 2$. Define a symmetric bilinear form $L$ on $\Sigma$ by
$$L(X, Y)= f'' g(X, Y) + A(X, Y) \nabla_\nu f$$
for any tangent vector fields $X$ and $Y$ on $\Sigma$. Then
\begin{align*}
\hess_{\scriptscriptstyle\Sigma} f (e_i, e_i) \leq L(e_i, e_i) = f'' + \lambda_i \nabla_\nu f,
\end{align*}
where each $\lambda_i$ is the principal curvature of $\Sigma$ for $i=1,2$. Since $\lambda_1+\lambda_2=2H$ and $\lambda_1^2 +\lambda_2^2=|A|^2$, using (\ref{eqn:squred norm}) gives
\begin{align} \label{eqn:gauss}
\lambda_1 \lambda_2 &= 2H^2 - \frac{1}{2} |A|^2 \nonumber \\
&=H^2- \frac{1}{2} |\Phi|^2.
\end{align}
Thus
\begin{align} \label{ineq: determinant}
\begin{array}{rl}
\det L & =\left(f''+\lambda_1 \nabla_\nu f\right)\left(f''+\lambda_2 \nabla_\nu f\right) \\
&= f''^2 +2Hf'' \nabla_\nu f + \lambda_1 \lambda_2 \left(\nabla_\nu f\right)^2  \\
&= \left( f'' +H \nabla_\nu f\right)^2 - \frac{1}{2} |\Phi|^2 (\nabla_\nu f)^2  \\
&\geq 0,
\end{array}
\end{align}
where we used the assumption (\ref{ineq:Assumption}) in the above inequality. We note that $f''<0$ if $c>0$ by definition of the function $f$. In order to prove $\hess_{\scriptscriptstyle \Sigma} f (e_i, e_i) \leq 0$ for $i=1, 2$, it suffices to show that $\tr L \leq 0$, since we have shown that $\det L \geq 0$ in (\ref{ineq: determinant}). Define a function $u$ on $\Sigma$ as follows:
$$u=f''+ H \nabla_\nu f =\frac{1}{2} \tr L.$$

Suppose that $H=0$. Because $c>0$, we see that $u=f'' <0$, which shows that $\tr L < 0$. Thus we may assume that $H\neq 0$. Now we claim that if $u(p)=0$ at some point $p \in \Sigma$,  then $|\Phi|(p) =0$, i.e., $p$ is an umbilical point. To see this, observe that, at the point $p\in \Sigma$ such that $u(p)=0$,
$$\nabla_\nu f (p) \neq 0$$
since $f''<0$ and $H\neq 0.$
Thus the assumption (\ref{ineq:Assumption}) implies that
$$|\Phi|^2(p) (\nabla_\nu f (p))^2 \leq u(p)^2=0,$$
which shows that $|\Phi|(p) = 0$. We remark that such umbilical points are isolated unless $\Sigma$ is totally umbilical by our assumption.

We claim that $\tr L$ has one sign, i.e., the function $u$ cannot change sign. To see this, suppose that $u$ changes sign. Then there exist two points $q_1, q_2 \in \Sigma$ such that $u(q_1)>0$ and $u(q_2)<0$. Choose any continuous curve $\gamma: [a, b] \rightarrow \Sigma$ with $\gamma(a)=q_1$ and $\gamma(b)=q_2$. By continuity, there exists a point $q=\gamma(t_0)$ for some $t_0 \in (a,b)$ such that $u(q)=0$ and $u$ changes sign near $q$. For simplicity, we may assume that

\begin{displaymath}
\left\{
\begin{array}{lll}
u(\gamma(t))>0  &~~\text{if } t\in [t_0-\varepsilon, t_0),\\
u(\gamma(t))=0  &~~\text{if } t= t_0,\\
u(\gamma(t))<0   &~~\text{if } t\in (t_0, t_0+\varepsilon]
\end{array}
\right.
\end{displaymath}
for some constant $\varepsilon>0$. Since $q$ is an isolated umbilical point of $\Sigma$, we can choose a geodesic ball $B_{r_0} (q)$ centered at $q$ with radius $r_0$ such that $q$ is the only umbilical point in $B_{r_0} (q)$. Choose a sufficiently small geodesic ball $B_{r_1} (q)$ for $r_1<r_0$ such that neither $\gamma(t_0-\varepsilon)$ nor $\gamma(t_0+\varepsilon)$ is contained in $B_{r_1} (q)$. Obviously, there is no umbilical point in the annulus $B_{r_0} (q)\setminus B_{r_1} (q)$. However, if we choose any two points $\gamma(t_0-\delta), \gamma(t_0+\delta) \in B_{r_0} (q)\setminus B_{r_1} (q)$ for some $0<\delta<\varepsilon$, then any continuous curve in $B_{r_0} (q)\setminus B_{r_1} (q)$ connecting the points $\gamma(t_0-\delta)$ and $\gamma(t_0+\delta)$ contains an umbilical pont in $B_{r_0} (q)\setminus B_{r_1} (q)$ by continuity, which is a contradiction. Therefore the function $u$ cannot change sign, i.e., $\tr L$ has one sign.

In order to finish the proof of the claim that $\tr L \leq0$, it suffices to prove that there exists a point at which the function $u$ is negative, since $u$ has one sign. Choose a point $p\in\Sigma$ such that the Gaussian curvature $K_\Sigma (p)<0$. Suppose $H\nabla_\nu f (p) \leq 0$. Then it follows that $u(p)<0$. Suppose $H\nabla_\nu f (p) >0$. Then, by (\ref{eqn:gauss}) and the assumption (\ref{ineq:Assumption}), the following inequalities hold:
\begin{displaymath}
\begin{array}{rcl}
|H| &<& \frac{1}{\sqrt{2}} |\Phi|(p).\\
\frac{1}{\sqrt{2}} |\Phi|(p) |\nabla_\nu f (p)| &\leq& \left|f''(p) +H  \nabla_\nu f(p) \right|.
\end{array}
\end{displaymath}
Combining these two inequalities,
$$0<H \nabla_\nu f(p) \leq \left|f''(p) +H  \nabla_\nu f (p)\right|.$$
Since $f''(p)<0$, we get
$$H \nabla_\nu f(p) < \left|f''(p)\right|,$$
which shows $u(p)<0$. Therefore we conclude that
$$\hess_{\scriptscriptstyle \Sigma} f (e_i, e_i) \leq L(e_i, e_i) \leq 0$$
for $i=1, 2$. In the same manner, one can prove $\hess_{\scriptscriptstyle \Sigma} f (e_i, e_i) \geq 0$ in case where $c\leq0$.\\

Before we move to the next step, let us define
$$\mathcal{A}=\{x\in \Sigma: r(x) =\min_{\scriptscriptstyle \Sigma} r\}.$$
Observe that
\begin{itemize}
  \item $\mathcal{A}\neq \emptyset$,
  \item $\mathcal{A} \cap \partial C = \emptyset$,
\end{itemize}
 since $\Sigma$ is a compact surface with free boundary in a strictly convex domain $C$.\\

\noindent Step IV: We claim that \emph{$\mathcal{A}$ is a totally convex subset of $\Sigma$}, i.e., any geodesic in $\Sigma$ connecting two points in $\mathcal{A}$ is contained in $\mathcal{A}$. To see this, let $\gamma: [a,b]\rightarrow \Sigma$ be a geodesic joining two points $q_1, q_2 \in \mathcal{A}$ such that $\gamma(a)=q_1$ and $\gamma(b)=q_2$. By Step III,

\begin{equation*}
\begin{cases}
0 \geq \frac{d^2}{dt^2} \left(f \circ \gamma \right)= \hess_{\scriptscriptstyle\Sigma} f(\frac{d\gamma}{dt}, \frac{d\gamma}{dt}) &~~{\text{if }} c > 0,\\
0 \leq \frac{d^2}{dt^2} \left(f \circ \gamma \right) &~~{\text{if }} c \leq 0.
\end{cases}
\end{equation*}
Since $q_1, q_2 \in \mathcal{A}$, we see that
$$\frac{d}{dt}\left( f\circ\gamma\right) (a)= \frac{d}{dt}\left( f\circ \gamma \right) (b)=0,$$
which implies that $f$ is constant on $\gamma$ by the maximum principle. Therefore $\gamma \subset \mathcal{A}$. Moreover, since $\partial \Sigma$ is strictly convex in $\Sigma$ by Step I, there exists a geodesic in $\Sigma$ connecting any two points $p_1,p_2 \in \Sigma$. This shows that the totally convex subset $\mathcal{A}$ is connected.\\

\noindent Step V: We claim that \emph{$\Sigma$ is homeomorphic to either a disk or an annulus}. To see this, we divide into two cases:\\

\noindent Case 1: $\mathcal{A}$ consists of a single point $q\in \Sigma\setminus \partial \Sigma$.\\
Case 2: $\mathcal{A}$ contains more than one point. \\

For Case 1, suppose that $[\alpha]$ is any non-trivial homotopy class in $\pi_1(\Sigma, q)$. Since $\partial \Sigma$ is strictly convex, there exists a closed geodesic $\gamma \subset \Sigma\setminus \partial \Sigma$ passing through $q$ which is homotopic to $\alpha$. Because $\mathcal{A}$ is totally convex by Step IV, we see that $\gamma \subset \mathcal{A} = \{q\}$, which is a contradiction. This implies that  $\pi_1(\Sigma, q)$ is trivial. Therefore $\Sigma$ is homeomorphic to a disk.

For Case 2, we may assume that $\Sigma$ is not homeomorphic to a disk. Choose a point $q\in \mathcal{A}$ and a non-trivial homotopy class $[\alpha] \in \pi_1(\Sigma, q)$. As before, since $\partial \Sigma$ is strictly convex, there exists a closed geodesic $\gamma \subset \Sigma\setminus \partial \Sigma$ passing through $q$ which is homotopic to $\alpha$. Note that $\gamma \subset \mathcal{A}$ by Step IV. For any $q_1\in \gamma$ which is different from $q$, we are going to show that any minimizing geodesic joining $q_1$ to an arbitrary point $q_2$ on $\gamma$ is contained in $\mathcal{A}$. Suppose that a minimizing geodesic $\tilde{\gamma}$ connecting $q_1$ and $q_2$ is not lying on $\gamma$. Since  $\mathcal{A}$ is totally convex and $\gamma \subset \mathcal{A}$, we see that $\tilde{\gamma} \subset \mathcal{A}$ and there exists a nonempty open subset $\mathcal{U}$ of $\mathcal{A}$. Thus, for any geodesic $\beta(t)$ in $\mathcal{U}$,
\begin{displaymath}
0= \frac{d^2}{dt^2} \left(f \circ \beta \right) =\hess_{\scriptscriptstyle\Sigma} f\left(\frac{d \beta}{dt}, \frac{d \beta}{dt}\right) =L\left(\frac{d \beta}{dt}, \frac{d \beta}{dt}\right)
\end{displaymath}
by Step III. This implies that
\begin{equation*}
\begin{cases}
L=0,\\
f'' +\lambda_1 \nabla_\nu f = f''+\lambda_2 \nabla_\nu f = 0
\end{cases}
\end{equation*}
in $\mathcal{U}$. Since $\nabla_\nu f \neq 0$, we have
$$\lambda_1=\lambda_2=-\frac{f''}{\nabla_\nu f}$$
in $\mathcal{U}$. Thus the open subset $\mathcal{U}$ is totally umbilical, which shows that $\Sigma$ must be totally umbilical by our assumption. However this is impossible in our situation. Therefore $\mathcal{A}$ has to be equal to the unique closed geodesic $\gamma$. It follows that $\Sigma$ is homeomorphic to an annulus since $[\alpha]$ was chosen to be arbitrary.

\end{proof}

\section{A gap theorem for free boundary cmc-$H$ surfaces in space forms}
In the previous section, we investigated the topology of free boundary cmc-$H$ surfaces in a strictly convex domain in a $3$-dimensional Riemannian manifold under some geometric condition on the length of the traceless second fundamental form. Restricting our attention to free boundary cmc-$H$ surfaces in a geodesic ball in space forms, we are able to obtain a gap theorem for such surfaces, which can be regarded as a characterization of spherical caps and Delaunay surfaces. To do this, we need some facts about rotation hypersurfaces in space forms.

Let $\Sigma$ be a hypersurface in an $(n+1)$-dimensional space form $\overline{M}^{n+1}(c)$ of constant curvature $c$. Choose any local orthonormal frame $\{e_1, \dots, e_{n}, e_{n+1}\}$ in $\overline{M}^{n+1}(c)$ such that, when restricted to $\Sigma$, the vectors
$e_1, \dots, e_{n-1}, e_{n}$ are tangent to $\Sigma$ and, consequently, $e_{n+1}$ is normal to $\Sigma$. We shall make use of the following convention on the range of the indices:
$$1 \leq A, B, C, \cdots \leq n+1, \ \ \ \ 1 \leq i, j, k, \cdots \leq n.$$
Let $\{\omega_A\}$ and $\{\omega_{AB}\}$ be the field of dual frames to $\{e_A\}$ and the connection forms on $\overline{M}^{n+1}(c)$, respectively. Then the structure equations of $\overline{M}^{n+1}(c)$ are given by
\begin{align}
d \omega_A=-\sum \omega_{AB} \wedge \omega_B&, ~~~ \omega_{AB}+\omega_{BA}=0, \nonumber\\
d \omega_{AB} =-\sum \omega_{AC} \wedge \omega_{CB} +\Omega_{AB}&, ~~~\Omega_{AB}=\frac{1}{2} \sum K_{ABCD}~ \omega_C \wedge \omega_D, \label{eqn:str eqn}\\
K_{ABCD} = c(\delta_{AC} \delta_{BD} &- \delta_{AD} \delta_{BC}),  \nonumber
\end{align}
where we used the Einstein's summation convention. Restricting these forms to $\Sigma$,
$$\omega_{n+1}=0.$$
Since
\begin{align*}
0=d \omega_{n+1} = - \sum \omega_{n+1, i} \wedge \omega_{i},
\end{align*}
by Cartan's lemma we may write
\begin{align}
\omega_{n+1, i} = \sum h_{ij} \omega_j, ~~~ h_{ij}=h_{ji}. \label{eqn: coeff}
\end{align}

From now on, we assume that $\Sigma$ is a rotation cmc-$H$ hypersurface in $\overline{M}^{n+1}(c)$, i.e., $\Sigma$ is a rotation hypersurface with constant mean curvature $H$. Parametrize $\Sigma$ by $f(t_1,\dots, t_{n-1},s)$ such that the directions of the parameters $t_1, \dots, t_{n-1},s$ are all principal directions and the $s$-parameter curve is a meridian. Let $\lambda_1,\dots,\lambda_{n-1}$ be the principal curvatures along the coordinate curves $t_i$ and let $\lambda_n=\mu$ be the principal curvature along the coordinate curve $s$ (see \cite{dCD} for instance).  Then $\Sigma$ has at most two distinct principal curvatures $\lambda$ and $\mu$ such that
\begin{equation*}
\begin{cases}
\lambda=\lambda(s)=\lambda_1=\dots=\lambda_{n-1},\\
\mu=\mu(s)=\lambda_n.
\end{cases}
\end{equation*}
\noindent Indeed, if a rotation cmc-$H$ hypersurface in $\overline{M}^{n+1}(c)$ has no umbilical point, then $\Sigma$ has two distinct principal curvatures.

Let $\alpha(s)$ be a meridian of $\Sigma$ parametrized by arclength parameter $s$. Choose a local orthonormal frame $\{e_1, \dots, e_{n}, e_{n+1}\}$ in $\overline{M}^{n+1}(c)$ such that, when restricted to $\Sigma$, the vectors
$e_1, \dots, e_{n-1}, e_{n}=\frac{\partial}{\partial s}$ are tangent to $\Sigma$ such that
$$h_{ij}=\lambda_i \delta_{ij}.$$
Then (\ref{eqn: coeff}) becomes
\begin{align}
\omega_{n+1, i} = \sum h_{ij} \omega_j = \lambda_i \omega_i. \label{eqn: cartan}
\end{align}
We take exterior differentiation of (\ref{eqn: cartan}) and define $h_{ijk}$ by
\begin{displaymath}
\sum h_{ijk} \omega_k =dh_{ij}-\sum  h_{ik}\omega_{kj}-\sum  h_{kj}\omega_{ki}.
\end{displaymath}
Then it follows that
$$h_{ijk}=h_{ikj}.$$
Note that, for each $1 \leq i \leq n-1$,
\begin{align} \label{eqn:directional derivative}
\nabla^{\Sigma}_{e_i} \lambda =e_i \lambda=0, ~~~\nabla^{\Sigma}_{e_i} \mu =e_i \mu=0.
\end{align}
Let $\theta_{ij}:=(\lambda_i -\lambda_j) \omega_{ij} = \theta_{ji}$. Then we have
\begin{align} \label{eqn: hijk}
\sum h_{ijk} \omega_k =\delta_{ij} d\lambda_j-(\lambda_i -\lambda_j) \omega_{ij}=\delta_{ij} d\lambda_j -\theta_{ij}.
\end{align}
From (\ref{eqn: hijk}) and the fact that $\lambda_i=\lambda_j=\lambda$ and $e_i \lambda=0$ for $1 \leq i, j \leq n-1$, it follows that, for $1\leq i \neq j \leq n-1$ and $1 \leq k \leq n$,
\begin{itemize}
  \item $h_{ijk}=h_{iij}=h_{nnj}=0$,
  \item $h_{iin}=e_n \lambda$, ~~$h_{nnn}=e_n \mu$.
\end{itemize}
\noindent Thus
\begin{displaymath}
\theta_{in}=\delta_{in} d\mu -\sum_{k=1}^n h_{ink} \omega_k  =-\sum_{k=1}^{n-1} h_{ink}\omega_k -h_{inn} \omega_n=-(e_n \lambda) \omega_i
\end{displaymath}
for $1\leq i\leq n-1$. Further, suppose that $\Sigma$ has no umbilical point. Obviously, $\lambda \neq \mu$ and $\lambda \neq H$. Then (\ref{eqn: hijk}) and the structure equations show that
\begin{align}
\omega_{in}&=\frac{\theta_{in}}{\lambda-\mu}=-\frac{e_n \lambda}{\lambda-\mu} \omega_i =- \frac{e_n \lambda}{n(\lambda-H)} \omega_i, \label{eqn:omega_in a}\\
d \omega_n &=-\sum \omega_{n i} \wedge \omega_i=0. \label{eqn:omega_in b}
\end{align}
From (\ref{eqn:omega_in b}), we see that
\begin{displaymath}
\omega_n =ds.
\end{displaymath}

Define a function $w$ by
$$w=\left|\lambda - H\right|^{-\frac{1}{n}}=\left(a\left(\lambda - H\right)\right)^{-\frac{1}{n}},$$
where $a=\text{sgn}(\lambda-H) \neq 0$. The directional derivative of $w$ with respect to $e_i$ for $1 \leq i \leq n$ is given by
\begin{align} \label{eqn:directional derivative w}
\nabla^{\Sigma}_{e_i} w=\nabla^{\Sigma}_{e_i} (|\lambda -H|^{-\frac{1}{n}})= -\frac{1}{n} a w^{n+1} \nabla^{\Sigma}_{e_i} \lambda .
\end{align}
For $i=1,\dots,{n-1}$, (\ref{eqn:directional derivative}) implies
\begin{equation*}
\nabla^{\Sigma}_{e_i} w=0.
\end{equation*}
Moreover, from (\ref{eqn:omega_in a}) and (\ref{eqn:directional derivative w}), we get
\begin{displaymath}
\omega_{i n} =- a\left(\log |\lambda -H|^{\frac{1}{n}} \right)' \omega_i =  \left(\log w \right)' \omega_i,
\end{displaymath}
where the prime denotes the derivative with respect to $s$. For $1\leq i \leq n-1$, (\ref{eqn:str eqn}) becomes
\begin{eqnarray*}
d \omega_{in} &=& -\sum_{k=1}^n \omega_{ik} \wedge \omega_{kn} - \omega_{i, n+1} \wedge \omega_{n+1, n}+c\, \omega_i \wedge \omega_n\\
&=&-\left(\log w \right)'\sum_{k=1}^{n-1} \omega_{ik} \wedge\omega_k + \lambda \mu \omega_i \wedge ds +c\, \omega_i \wedge ds.
\end{eqnarray*}
On the other hand, a direct computation shows that
\begin{eqnarray*}
d \omega_{in} &=&\left(\log w \right)'' ds \wedge \omega_i +\left(\log w \right)' d\omega_i\\
&=&-\left(\log w \right)'' \omega_i \wedge ds -\left(\log w \right)' \sum_{k=1}^n \omega_{ik} \wedge \omega_k\\
&=&\left\{-\left(\log w \right)'' -\left(\left(\log w \right)'\right)^2 \right\}  \omega_i \wedge ds -\left(\log w \right)' \sum_{k=1}^{n-1} \omega_{ik} \wedge\omega_k.
\end{eqnarray*}
Therefore we obtain
\begin{equation} \label{eqn:govern}
-\left(\log w \right)'' -\left(\left(\log w \right)'\right)^2 =\lambda \mu +c.
\end{equation}
By a straightforward computation, we have
\begin{align*}
\lambda \mu &=\lambda \left(nH-(n-1)\lambda \right)\\
&=(1-n) (\lambda-H)^2 + (2-n)H(\lambda-H)+H^2\\
&=(1-n)w^{-2n}+a(2-n)Hw^{-n}+H^2,
\end{align*}
which enables us to reformulate (\ref{eqn:govern}) as follows:
\begin{equation*}
w''+w \left(c+H^2+a(2-n)Hw^{-n}+(1-n)w^{-2n} \right)=0.
\end{equation*}
In summary, we have the following proposition, which was already obtained when $H=0$  by Otsuki \cite{Otsuki} and when $c\geq 0$ by Wei \cite{Wei}.
\begin{prop}[\cite{Otsuki, Wei}]
Let $\overline{M}^{n+1}(c)$ be an $(n+1)$-dimensional space form of constant curvature $c$. Let $\Sigma$ be a rotation cmc-$H$ hypersurface in $\overline{M}^{n+1}(c)$. For any non-umbilical point $p \in \Sigma$, we have the following ordinary differential equation in a neighborhood of $p$:
\begin{equation} \label{eqn:govern2}
w''+w \left(c+H^2+a(2-n)Hw^{-n}+(1-n)w^{-2n} \right)=0,
\end{equation}
where $w=\left|\lambda - H\right|^{-\frac{1}{n}}$ and $a=\text{sgn}(\lambda-H)$.
\end{prop}

In \cite{Pyo}, Pyo proved that if $\Sigma$ is an immersed minimal surface in $\RR^3$ which meets a plane along a circle with constant angle, then $\Sigma$ is part of catenoid. Before stating a gap theorem for free boundary cmc-$H$ surfaces, we show the existence and uniqueness of cmc-$H$ surfaces in space forms meeting a totally geodesic plane orthogonally along a geodesic circle, which may be of interest in its own right. More precisely, we prove the following, which can be regarded as an extension of Pyo's result \cite{Pyo}.

\begin{lem} \label{lem:rotation}
Let $\overline{M}^3(c)$ be a $3$-dimensional space form of constant curvature $c$. Let $P$ be a totally geodesic plane in $\overline{M}^3(c)$ and $\gamma$ be a geodesic circle lying on $P$. Then, for given a constant $H$, there exists a unique cmc-$H$ surface $\Sigma \subset \overline{M}^3(c)$ which meets $P$ orthogonally along $\gamma$. In particular, such $\Sigma$ is a rotation cmc-$H$ surface, i.e., a Delaunay surface.
\end{lem}

\begin{proof}

(Uniqueness) We recall that any cmc-$H$ surface in space forms is locally the graph of an analytic function over the tangent space at some point. Indeed, this analytic function is the solution of the mean curvature equation, which is a second-order elliptic partial differential equation. Let $\Sigma_1$ and $\Sigma_2$ be cmc-$H$ surfaces that meet tangentially along the given geodesic circle $\gamma$. Then the classical maximum principle guarantees that $\Sigma_1 =\Sigma_2$.\\

\noindent (Existence)  Let  $\gamma \subset P$ be a geodesic circle of radius $r_0$ centered at a point $p\in P$. Note that any cmc surface has an isolated umbilical point unless it is totally umbilical. First, suppose that $\Sigma$ is totally umbilical. In this case, $\Sigma$ is a geodesic sphere in $\overline{M}^3(c)$. Choose $\Sigma = \partial B_{r_0}(p)$, where $B_{r_0}(p)$ denotes the geodesic ball of radius $r_0$ centered at $p$. Then $\Sigma$ meets $P$ orthogonally along $\gamma$. Moreover, $\gamma$ is a line of curvature of $\Sigma$ by Joachimsthal's theorem \cite{Spivak}. Note that the principal curvature $\lambda$ along $\gamma$ is constant and the mean curvature $H$ of $\Sigma$ satisfies $H=\lambda$.

Now suppose $\Sigma$ has an isolated umbilical point. We will show the existence of a rotation cmc surface $\Sigma$ in $\overline{M}^3(c)$ meeting $P$ orthogonally along $\gamma$ by solving an ordinary differential equation in terms of the principal curvature $\lambda$ of $\Sigma$. To do this, let us find the ordinary differential equation which is satisfied on such $\Sigma$ meeting $P$ orthogonally along $\gamma$. Since umbilical points on $\Sigma$ are isolated, we can choose a non-umbilical point $q$ on $\gamma \subset \Sigma$.  Then it follows from (\ref{eqn:govern2}) that in a neighborhood of $q$,
\begin{equation*}
w'' + w\left(c+H^2 -w^{-4}\right)=0,
\end{equation*}
where $w=|\lambda-H|^{-\frac{1}{2}}$. Equivalently,
\begin{equation} \label{eqn:w1}
\left(w'\right)^2 + w^2 \left(c+ \left(H +w^{-2} \right)^2 \right)=C,
\end{equation}
for some constant $C$. Here the differentiation is taken with respect to the arclength parameter $s$ of a meridian $\alpha(s)$. We may assume that $\alpha(0)= q \in \gamma$. Multiplying  both sides of (\ref{eqn:w1}) by $w^2$,
\begin{equation*}
w^2 \left(w' \right)^2 + \left(c+H^2 \right)w^4 - \left(C-2H \right)w^2+1=0.
\end{equation*}
By substituting $w^2$ with $u$, the above equation can be rewritten by the first-order ordinary differential equation
\begin{equation} \label{eqn:u1}
\frac{1}{4} \left(u' \right)^2 + \left(c+H^2 \right)u^2 - \left(C-2H \right)u+1=0,
\end{equation}
which gives
\begin{equation} \label{eqn:u2}
\displaystyle\frac{du}{2\sqrt{- \left(c+H^2 \right)u^2 +\left(C-2H \right)u - 1}} =ds.
\end{equation}
Since $u=w^2$, we have the following initial condition:
\begin{align} \label{eqn:initial1}
u_0:=u|_{s=0}=w^2|_{s=0}=|\lambda-H|^{-1} |_{s=0} > 0.
\end{align}
Note that $\lambda' \equiv 0$ on $\gamma$, since $\Sigma$ meets $P$ orthogonally along $\gamma$. Thus by using (\ref{eqn:directional derivative w}), we get another initial condition:
\begin{align} \label{eqn:initial2}
u'|_{s=0}=2ww'|_{s=0}=-aw^4\lambda' |_{s=0} = 0.
\end{align}
For $s=0$, substituting (\ref{eqn:initial1}) and (\ref{eqn:initial2}) into (\ref{eqn:u1}) yields
\begin{align} \label{eqn:algebraic}
\left(c+H^2 \right)u_0^2 - \left(C-2H \right)u_0+1=0.
\end{align}
We remark that the constant $C$ is completely determined by (\ref{eqn:algebraic}). On the the hand, if we consider (\ref{eqn:algebraic}) as a second-order algebraic equation in $u_0$, then the discriminant of (\ref{eqn:algebraic}) satisfies that
\begin{align} \label{eqn:discriminant}
\left(C-2H \right)^2-4\left(c+H^2 \right) =C^2-4HC-4c \geq 0.
\end{align}
Note that
\begin{displaymath}
- \left(c+H^2 \right)u^2 +\left(C-2H \right)u - 1 = -\left(c+H^2\right) \left( u -\displaystyle\frac{C-2H}{2(c+H^2)}      \right)^2 +\displaystyle\frac{(C-2H)^2}{4(c+H^2)} -1.
\end{displaymath}
In order to solve (\ref{eqn:u2}), we divide into the following three cases:\\

\noindent Case 1: $c+H^2>0$,\\
Case 2:  $c+H^2<0$, \\
Case 3: $c+H^2=0$. \\

For Case 1, from (\ref{eqn:discriminant}), it follows that
$C^2-4HC-4c>0$ and (\ref{eqn:u2}) becomes
\begin{equation*}
\displaystyle\frac{1}{2\sqrt{c+H^2}}\displaystyle\frac{du}{\sqrt{\displaystyle\frac{C^2-4HC-4c}{4(c+H^2)^2}-\left( u-\displaystyle\frac{C-2H}{2(c+H^2)}\right)^2}} =ds,
\end{equation*}
which gives
\begin{equation} \label{eqn:sol1}
u(s)= \displaystyle\frac{C-2H}{2(c+H^2)}+ {\displaystyle\frac{\sqrt{C^2-4HC-4c}}{2(c+H^2)} \sin (2\sqrt{c+H^2} s -D}),
\end{equation}
where $D$ is the constant satisfying that $\cos D=0$ by (\ref{eqn:initial2}).

For Case 2, using (\ref{eqn:discriminant}) again, we see that
$C^2-4HC-4c>0$. Moreover, by using (\ref{eqn:u2}) and (\ref{eqn:initial2}),
\begin{equation} \label{eqn:sol2}
u(s)=\displaystyle\frac{C-2H}{2(c+H^2)}+ {\frac{\sqrt{C^2-4HC-4c}}{2(c+H^2)} \cosh (2\sqrt{-(c+H^2)} s}).
\end{equation}

For Case 3, (\ref{eqn:u2}) becomes
\begin{equation*}
\displaystyle\frac{du}{2\sqrt{(C-2H)u - 1}} =ds.
\end{equation*}
Note that if $c+H^2=0$, then $C-2H =\frac{1}{u_0} > 0$ by (\ref{eqn:algebraic}). Then the solution is given by
\begin{equation} \label{eqn:sol3}
u(s)= \displaystyle\frac{1}{C-2H}+ \left(\sqrt{C-2H} s \right)^2.
\end{equation}

So far we have found the explicit solution $u(s)$ in any case, which shows that the principal curvatures $\lambda$ and $\mu$ are given explicitly. Since a rotation surface is determined by its principal curvatures (see \cite{dCD} for instance), we obtain a rotation cmc-$H$ surface $\Sigma$ with the given principal curvatures $\lambda$ and $\mu$. Moreover, the surface $\Sigma$ meets a totally geodesic plane $P$ orthogonally along the geodesic circle $\gamma$. This completes the proof.
\end{proof}

It seems interesting to prove the existence of cmc-$H$ surfaces in space forms meeting a totally geodesic plane along a geodesic circle with a constant angle which is not necessarily $\frac{\pi}{2}$. Before stating our theorem, we recall some basic facts concerning the Hopf differential.

Let $\Sigma$ be an immersed surface in a $3$-dimensional Riemannian manifold $M^3$. Suppose $(u, v)$ is a local isothermal coordinate system on $\Sigma$ and the second fundamental form $A$ is given by ${\cal L} du^2+2{\cal M} du dv+ {\cal N} dv^2$.  The Hopf differential $Q$ of $\Sigma$ is the $(2,0)$-part of the complexified second fundamental form written as
\begin{displaymath}
Q=\frac{1}{4} \left( ({\cal{L}} -{\cal{N}}) -2i {\cal{M}} \right) dz^2,
\end{displaymath}
where $z=u+iv$. It is well-known that if the Hopf differential $Q$ is holomorphic in $\Sigma$, then the zero set of $Q$ is isolated unless $Q \equiv 0$. We remark that a point $p \in \Sigma$ is umbilical if $Q(p)=0$. We are now in position to state a gap theorem for free boundary cmc-$H$ surfaces in geodesic balls of space forms.

\begin{thm} \label{thm:main}
Let $\overline{M}^3(c)$ be a $3$-dimensional space form of constant curvature $c$. Let $\Sigma$ be a free boundary cmc-$H$ surface in the geodesic ball $B_R(p) \subset \overline{M}^3(c)$ of radius $R$ centered at $p$. When $c>0$, assume that $R<\frac{\pi}{2 \sqrt{c}}$.  For the function $f$ defined as in (\ref{eqn:definition}) and the distance $r$ measured from $p$, assume that every point in $\Sigma$ satisfies
\begin{align} \label{ineq:Assumption2}
\frac{1}{2} |\Phi|^2 \left(\nabla_\nu f \right)^2 \leq \left(f'' +H  \nabla_\nu f \right)^2,
\end{align}
where $\nabla_\nu f = \langle \nu, \nabla f \rangle$. Here $\nu$ and $\nabla$ denote the unit normal vector field on $\Sigma$ and the Levi-Civita connection on $\overline{M}^3(c)$, respectively.  Then $\Sigma$ is either a spherical cap or a Delaunay surface.
\end{thm}

\begin{proof}
Since the Hopf differential $Q$ on $\Sigma$ in $\overline{M}^3(c)$ is holomorphic, we see that umbilic points on $\Sigma$ is isolated unless $\Sigma$ is totally umbilical. Thus $\Sigma$ is homeomorphic to a disk or an annulus by Theorem \ref{thm:topology}. Suppose that $\Sigma$ is homeomorphic to a disk. Then it is well-known that $\Sigma$ is a spherical cap (see \cite{Nitsche, RS, Souam}). Suppose that $\Sigma$ is homeomorphic to an annulus. Define
$$\mathcal{A}=\{x\in \Sigma: r(x) =\min_{\scriptscriptstyle \Sigma} r\}.$$
Applying the same argument as in Step V of the proof of Theorem \ref{thm:topology}, we obtain that there exists a unique closed geodesic $\gamma$ in $\Sigma$ satisfying that $\gamma=\cal{A}$. This implies that $\Sigma$ meets tangentially the geodesic sphere $\partial B_{r_0}(p)$ of radius $r_0$ centered at $p$ along $\gamma$, where $r_0=\min_{\scriptscriptstyle \Sigma} r$. From the Joachimsthal's Theorem \cite{Spivak}, it follows that $\gamma$ is a line of curvature of $\Sigma$. Since $\gamma$ is both a geodesic and a line of curvature of $\Sigma$, we conclude that $\gamma$ is the geodesic circle lying on a totally geodesic plane $P$ and $\Sigma$ meets $P$ orthogonally along $\gamma$. By Lemma \ref{lem:rotation}, $\Sigma$ is a Delaunay surface in $\overline{M}^3(c)$.

\end{proof}

As direct consequences of Theorem \ref{thm:main}, we have the following corollaries.

\begin{cor}
Let $\overline{M}^3(c)$ be a $3$-dimensional space form of constant curvature $c$. Let $\Sigma$ be a free boundary cmc-$H$ surface in the geodesic ball $B_R(p) \subset \overline{M}^3(c)$ of radius $R$ centered at $p$. When $c>0$, assume that $R<\frac{\pi}{2 \sqrt{c}}$.  For the function $f$ defined as in (\ref{eqn:definition}) and the distance $r$ measured from $p$, assume that every point in $\Sigma$ satisfies
\begin{align*}
\frac{1}{2} |\Phi|^2 \left(\nabla_\nu f \right)^2 < \left(f'' +H  \nabla_\nu f \right)^2.
\end{align*}
Then $\Sigma$ is a spherical cap.
\end{cor}

\begin{proof}
By assumption and the argument as in Step III of the proof of Theorem \ref{thm:topology}, we have, for any tangent vector field $v$ on $\Sigma$,
\begin{equation*}
\begin{cases}
\hess_{\scriptscriptstyle \Sigma} f (v, v) < 0  &~~{\text{if }} c > 0,\\
\hess_{\scriptscriptstyle \Sigma} f (v, v) > 0 &~~{\text{if }}  c \leq 0.
\end{cases}
\end{equation*}
Consider the case where $c>0$. Then the function $f$ is strictly concave, which implies that the set $\mathcal{A}=\{x\in \Sigma: r(x) =\min_{\scriptscriptstyle \Sigma} r\}$ consists of a single point. The argument as in Step V of the proof of Theorem \ref{thm:topology} shows that $\Sigma$ is homeomorphic to a disk. Therefore, by Theorem \ref{thm:main}, we conclude that $\Sigma$ is a spherical cap. In the same manner, one can get the same conclusion when $c\leq 0$.
\end{proof}

\begin{cor}
Let $\overline{M}^3(c)$ be a $3$-dimensional space form of constant curvature $c$. Let $\Sigma$ be a free boundary cmc-$H$ surface in the geodesic ball $B_R(p) \subset \overline{M}^3(c)$ of radius $R$ centered at $p$. When $c>0$, assume that $R<\frac{\pi}{2 \sqrt{c}}$.  For the function $f$ defined as in (\ref{eqn:definition}) and the distance $r$ measured from $p$, assume that every point in $\Sigma$ satisfies
\begin{align} \label{ineq: condition}
\frac{1}{2} |\Phi|^2 \left(\nabla_\nu f \right)^2 \leq \left(f'' +H  \nabla_\nu f \right)^2.
\end{align}
If equality in (\ref{ineq: condition}) holds at some point $q \in \Sigma$, then $\Sigma$ is a Delaunay surface.
\end{cor}

\begin{proof}
Suppose that $\Sigma$ is a spherical cap. Then $\Phi \equiv 0$ since $\Sigma$ is totally umbilical. If $f'=0$, then
$$f''+ H \nabla_\nu f = f''+ H f' \nabla_\nu r\neq 0$$
by definition of the function $f(r)$. Thus the right hand side of (\ref{ineq: condition}) never vanishes, which implies that equality in (\ref{ineq: condition}) cannot hold at any point in $\Sigma$.  Hence we may assume that $f' \neq 0$. Then
\begin{align*}
\left(f''+ H \nabla_\nu f\right)^2=\left(f'\right)^2 \left(\frac{f''}{f'}+H\nabla_\nu r\right)^2 {\rm ~~~and ~~~ } \frac{f''}{f'}>0.
\end{align*}
Observe that
$$H \nabla_\nu r \geq 0$$
on $\Sigma$, since $\Sigma$ is a spherical cap with free boundary in the geodesic ball $B_R (p)$. Thus equality in (\ref{ineq: condition}) cannot occur at any point in $\Sigma$, which is a contradiction. Therefore $\Sigma$ is a Delaunay surface by Theorem \ref{thm:main}.
\end{proof}

\vskip 0.3cm
\noindent
{\bf Acknowledgment:} The first author was supported by Basic Science Research Program through the National Research Foundation of Korea(NRF) funded by the Ministry of Education (Grant Number: 2017R1D1A1B03036369). The second author was supported by the National Research Foundation of Korea (NRF-2021R1A2C1003365).


\vskip 1cm
\noindent Sung-Hong Min\\
Department of Mathematics\\
Chungnam National University\\
Daehak-Ro 99, Yuseong-Gu, Daejeon, 34134, Korea \\
{\tt e-mail:sunghong.min@cnu.ac.kr}\\

\bigskip
\noindent Keomkyo Seo\\
Department of Mathematics and Research Institute of Natural Science\\
Sookmyung Women's University\\
Cheongpa-ro 47-gil 100, Yongsan-ku, Seoul, 04310, Korea \\
{\tt E-mail:kseo@sookmyung.ac.kr}\\
URL: http://sites.google.com/site/keomkyo/


\begin{thebibliography}{123}
\bibitem{AS} F. J. Almgren Jr., L. Simon, {\em Existence of embedded solutions of Plateau's problem}, Ann. Scuola Norm. Sup. Pisa Cl. Sci. (4) {\bf 6} (1979), no. 3, 447-495.
\bibitem{AN} L. Ambrozio, I. Nunes, {\em A gap theorem for free boundary minimal surfaces in the three-ball}, Comm. Anal. Geom. {\bf 29} (2021) no. 2, 283-292.
\bibitem{BCP} E. Barbosa, M. P. Cavalcante, E. Pereira, {\em Gap results for free boundary CMC surfaces in the Euclidean three-ball}, arXiv:1908.09952 [math.DG] (2019).
\bibitem{BV} E. Barbosa, C. Viana, {\em A remark on a curvature gap for minimal surfaces in the ball}, Math. Z. {\bf 294} (2020), no. 1-2, 713–720.
\bibitem{BB} J. Bernstein, C. Breiner, {\em A variational characterization of the catenoid}, Calc. Var. Partial Differential Equations {\bf 49} (2014), no. 1-2, 215–232.
\bibitem{CFS} A. Carlotto, G. Franz, M. B. Schulz, {\em Free boundary minimal surfaces with connected boundary and arbitrary genus}, arXiv:2001.04920 [math.DG].
\bibitem{CMV} M. P. Cavalcante, A. Mendes, F. Vit\'{o}rio, {\em Vanishing theorems for the cohomology groups of free boundary submanifolds}, Ann. Global Anal. Geom. {\bf 56} (2019), no. 1, 137–146.
\bibitem{CDK} S. S. Chern, M. do Carmo, S. Kobayashi, {\em Minimal submanifolds of a sphere with second fundamental form of constant length}, 1970, Functional Analysis and Related Fields (Proc. Conf. for M. Stone, Univ. Chicago, Chicago, Ill., 1968) 59–75 Springer, New York.
\bibitem{Courant} R. Courant, {\em The existence of minimal surfaces of given topological structure under prescribed boundary conditions}, Acta Math. {\bf 72} (1940), 51–98.
\bibitem{Devyver} B. Devyver, {\em Index of the critical catenoid}, Geom. Dedicata {\bf 199} (2019), 355–371.
\bibitem{dCD} M. Do Carmo, M. Dajczer, {\em Rotation hypersurfaces in spaces of constant curvature}, Trans. Amer. Math. Soc., {\bf 277} (1983), {no. 2}, 685-709.
\bibitem{GW} R. E. Greene, H. Wu, {\em Function theory on manifolds which possess a pole}, Lecture Notes in Mathematics, {\bf 699}, Springer, Berlin, 1979.
\bibitem{GJ} M. Gr\"{u}ter, J.  Jost, {\em On embedded minimal disks in convex bodies}, Ann. Inst. H. Poincar\'{e} Anal. Non Lin\'{e}aire {\bf 3} (1986), no. 5, 345–390.
269--287.
\bibitem{FPZ} A. Folha, F. Pacard, T. Zolotareva, {\em Free boundary minimal surfaces in the unit $3$-ball}, Manuscripta Math. {\bf 154} (2017), no. 3-4, 359–409.
\bibitem{Fraser} A. Fraser, {\em On the free boundary variational problem for minimal disks}, Comm. Pure Appl. Math. {\bf 53} (2000), no. 8, 931–971.
\bibitem{FS2011} A. Fraser, R. Schoen, {\em The first Steklov eigenvalue, conformal geometry, and minimal surfaces}, Adv. Math. {\bf 226} (2011), no. 5, 4011–4030.
\bibitem{FS} A. Fraser, R. Schoen, {\em Uniqueness theorems for free boundary minimal disks in space forms}, Int. Math. Res. Not. {\bf 2015}, no. 17, 8268–8274.
\bibitem{FS2016} A. Fraser, R. Schoen, {\em Sharp eigenvalue bounds and minimal surfaces in the ball}, Invent. Math. {\bf 203} (2016), no. 3, 823–890.
\bibitem{FS2020} A. Fraser, R. Schoen, {\em Some results on higher eigenvalue optimization}, Calc. Var. Partial Differential Equations {\bf 59} (2020), no. 5, Paper No. 151, 22 pp.
\bibitem{Hildebrandt} S. Hildebrandt, {\em Free boundary problems for minimal surfaces and related questions}, Frontiers of the mathematical sciences: 1985 (New York, 1985). Comm. Pure Appl. Math. {\bf 39} (1986), no. S, suppl., S111–S138.
\bibitem{KL} N. Kapouleas, M. Li, {\em Free boundary minimal surfaces in the unit three-ball via desingularization of the critical catenoid and the equatorial disk}, arXiv:1709.08556 [math.DG].
\bibitem{KW} N. Kapouleas, D. Wiygul, {\em Free-boundary minimal surfaces with connected boundary in the $3$-ball by tripling the equatorial disc}, arXiv:1711.00818v2 [math.DG].
\bibitem{Ketover} D. Ketover, {\em Free boundary minimal surfaces of unbounded genus}, 	arXiv:1612.08691 [math.DG].
\bibitem{KM} R. Kusner, P. McGrath, {\em On Free boundary minimal annuli embedded in the unit ball}, arXiv:2011.06884 [math.DG].
\bibitem{Lawson} H. B. Lawson Jr., {\em Local rigidity theorems for minimal hypersurfaces}, Ann. of Math. (2) {\bf 89} (1969), 187–197.
\bibitem{Lewy} H. Lewy, {\em On mimimal surfaces with partially free boundary}, Comm. Pure Appl. Math. {\bf 4} (1951), 1–13.
\bibitem{Li} M. Li, {\em Free boundary minimal surfaces in the unit ball : recent advances and open questions}, arXiv:1907.05053v3 [math.DG], to appear in Proceedings of the first annual meeting of the ICCM.
\bibitem{LX} H. Li, C. Xiong, {\em A gap theorem for free boundary minimal surfaces in geodesic balls of hyperbolic space and hemisphere}, J. Geom. Anal., {\bf 28} (2018), no. 4, 3171-3182.
\bibitem{LZ} M. Li, X. Zhou, {\em Min-max theory for free boundary minimal hypersurfaces I - regularity theory}, arXiv:1611.02612v3 [math.DG].
\bibitem{MNS} D. Maximo, I. Nunes, G. Smith, {\em Free boundary minimal annuli in convex three-manifolds}, J. Differential Geom. {\bf 106} (2017), no. 1, 139–186.
\bibitem{McGrath} P. McGrath, {\em A characterization of the critical catenoid}, Indiana Univ. Math. J. {\bf 67} (2018), no. 2, 889–897.
\bibitem{MPR} W. H. Meeks III, J. P\'{e}rez, A. Ros, {\em Stable constant mean curvature surfaces}, Handbook of geometric analysis. no. 1, 301–380, Adv. Lect. Math. (ALM), {\bf 7}, Int. Press, Somerville, MA, 2008.
\bibitem{Nitsche} J. C. C. Nitsche, {\em Stationary partitioning of convex bodies.}, Arch. Rational Mech. Anal. {\bf 89} (1985), no. 1, 1–19.
\bibitem{Otsuki} T. Otsuki, {\em Minimal hypersurfaces in a {R}iemannian manifold of constant curvature}, Amer. J. Math., {\bf 92} (1970), 145-173.
\bibitem{Petersen} P. Petersen, {\em Riemannian geometry}, Third edition. Graduate Texts in Mathematics, {\bf 171}, Springer, Cham, 2016.
\bibitem{Pyo} J. Pyo, {\em Minimal annuli with constant contact angle along the planar boundaries}, Geom. Dedicata, {\bf 146} (2010), 159-164.
\bibitem{RS} A. Ros, R. Souam, {\em On stability of capillary surfaces in a ball}, Pacific J. Math., {\bf 178} (1997), 345-361.
\bibitem{SY} R. Schoen, S.-T. Yau, {\em Lectures on differential geometry}, Conference Proceedings and Lecture Notes in Geometry and Topology, I. International Press, Cambridge, MA, 1994.
\bibitem{Simons} J. Simons, {\em Minimal varieties in riemannian manifolds}, Ann. of Math. (2) {\bf 88} (1968), 62–105.
\bibitem{SZ} G. Smith, D. Zhou, {\em The Morse index of the critical catenoid}, Geom. Dedicata {\bf 201} (2019), 13–19.

\bibitem{Souam} R. Souam, {\em On stability of stationary hypersurfaces for the partitioning problem for balls in space forms}, Math. Z. {\bf 224} (1997), no. 2, 195–208.
\bibitem{Spivak} M. Spivak, {\em A comprehensive introduction to differential geometry. Vol. III}, Second edition, Publish or Perish, Inc., Wilmington, Del., 1979.
\bibitem{Struwe} M. Struwe, {\em On a free boundary problem for minimal surfaces},  Invent. Math. {\bf 75} (1984), no. 3, 547-560.
\bibitem{Tran} H. Tran, {\em Index characterization for free boundary minimal surfaces}, Comm. Anal. Geom. {\bf 28} (2020), no. 1, 189–222.
\bibitem{Wei} G. Wei, {\em Complete hypersurfaces with constant mean curvature in a unit sphere}, Monatsh. Math., {\bf 149} (2006), {no. 3}, 251-258.
\bibitem{WW} G. Wheeler, V.-M. Wheeler, {\em Minimal hypersurfaces in the ball with free boundary}, Differential Geom. Appl. {\bf 62} (2019), 120–127.
\end{thebibliography}
\end{document}